\documentclass[a4paper]{article}
\usepackage{amssymb}
\usepackage{amsthm}
\usepackage{tikz}
\usepackage{mathtools}
\usepackage[title,toc,page,header]{appendix}
\usepackage{titletoc}

\usepackage[mathscr]{euscript} 
\DeclareFontFamily{OT1}{pzc}{}
\DeclareFontShape{OT1}{pzc}{m}{it}{<-> s * [1.10] pzcmi7t}{}
\DeclareMathAlphabet{\mathpzc}{OT1}{pzc}{m}{it} 

\dottedcontents{section}[2.5em]{\bfseries}{2.5em}{1pc}

\usetikzlibrary{calc}
\usetikzlibrary{matrix,arrows,decorations}
\newtheorem{theorem}{Theorem}[section]
\newtheorem{lemma}[theorem]{Lemma}

\newtheorem{corollary}[theorem]{Corollary}
\newtheorem{remark}[theorem]{Remark}
\newtheorem{definition}[theorem]{Definition}

\newtheorem*{conjecture*}{Conjecture}
\newtheorem*{theorem*}{Theorem}
\newtheorem*{remark*}{Remark}
\usepackage{hyperref}

\title{Fat realization and Segal's classifying space}

\author{Yi-Sheng Wang}
 
\begin{document}
\maketitle
\begin{abstract}
In this paper, we give a new proof of a well-known theorem due to tom Dieck that the fat realization and Segal's classifying space of an internal category in the category of topological spaces are homotopy equivalent.  
\end{abstract}

\tableofcontents 
\section{Introduction} 
Motivated by bundle theory, foliation theory, and delooping theory, classifying spaces of topological groups and groupoids were intensively studied during the 60s, 70s and 80s. Since then, many different constructions of classifying spaces of topological groups and groupoids have been introduced, for example, the Milnor construction, the Segal construction, fat realization, geometric realization and so on \cite{Mil2}, \cite{Mil3}, \cite{DL}, \cite{Sta1}, \cite{Sta2}, \cite{Sta3}, \cite{Se1}, \cite{Ha}, \cite{Bo}, \cite{Se3}, \cite{tD}. Some of them have even been generalized to any internal categories in $\mathpzc{Top}$, the category of topological spaces. For topological groups, most of the constructions give rise to homotopy equivalent spaces. However, for general internal categories in $\mathpzc{Top}$, the relation between them is not always clear. In this paper, we shall focus on the comparison between the Segal construction and fat realization of internal categories in $\mathpzc{Top}$.

\subsection*{Main results}
tom Dieck's theorem \cite[Proposition $2$]{tD} asserts that, given any simplicial space $X_{\cdot}$, the projection
\begin{equation}\label{Intro:tomDieckthm}
\pi:\vert\vert X\times S_{\cdot}\vert\vert\rightarrow \vert\vert X_{\cdot}\vert\vert
\end{equation}
is a homotopy equivalence, where $S_{\cdot}$ is the semi-simplicial set given by strictly increasing sequences of natural numbers. The idea of his proof is to construct a map $\vert\vert X_{\cdot}\vert\vert\rightarrow\vert\vert X\times S_{\cdot}\vert\vert$ and show it is a homotopy inverse to the projection $\pi$. However, to find a well-defined map from $\vert\vert X_{\cdot}\vert\vert$ to $\vert\vert X_{\cdot}\times S_{\cdot}\vert\vert$ and show it is a homotopy inverse to $\pi$ turn out to be quite difficult and complicated\footnote{The homotopy inverse $\rho$ given in \cite{tD} appears not to be well-defined. See the appendix for more details.}. To get around this inconvenience, we employ Quillen's theorem $A$. Our approach is more conceptual, but it works only when $X_{k}$ has the homotopy type of a $\operatorname{CW}$-complex and $X_{\cdot}=\operatorname{Ner}_{\cdot}\mathcal{C}$, where $\operatorname{Ner}_{\cdot}\mathcal{C}$ is the nerve of an internal category in $\mathpzc{Top}$.
  
\begin{theorem}\label{Thm1}
Given $\mathcal{C}$ an internal category in $\mathpzc{Top}$ such that its nerve $\operatorname{Ner}_{\cdot}\mathcal{C}$ has the homotopy type of a $\operatorname{CW}$-complex at each degree, then the canonical projection
\[\pi:\vert\vert \operatorname{Ner}_{\cdot}\mathcal{C} \times S_{\cdot}\vert\vert\rightarrow \vert\vert\operatorname{Ner}_{\cdot}\mathcal{C}\vert\vert\]
is a homotopy equivalence. 
\end{theorem}
Combing Theorem \ref{Thm1} with the fact that Segal's classifying space $\vert\operatorname{Ner}_{\cdot}\mathcal{C}^{\mathbb{N}}\vert$ is homeomorphic to the space $\vert\vert\operatorname{Ner}_{\cdot}\mathcal{C}\times S_{\cdot}\vert\vert$, where $\mathcal{C}^{\mathbb{N}}$ is Segal's unraveled category of $\mathcal{C}$ over the natural numbers $\mathbb{N}$, one can easily deduce the following useful corollary.
\begin{corollary}
If, in addition to the assumptions in Theorem \ref{Thm1}, the simplicial space $\operatorname{Ner}_{\cdot}\mathcal{C}$ is proper, then the forgetful functor 
\[\mathcal{C}^{\mathbb{N}}\rightarrow \mathcal{C}\]
induces a homotopy equivalence
\[\vert\operatorname{Ner}_{\cdot}\mathcal{C}^{\mathbb{N}}\vert\rightarrow \vert\operatorname{Ner}_{\cdot}\mathcal{C}\vert.\]
\end{corollary}
The space $\vert\operatorname{Ner}_{\cdot}\mathcal{C}^{\mathbb{N}}\vert$ is more natural from the point of view of bundle theory, whereas, category-theoretically, the geometric realization $\vert\operatorname{Ner}_{\cdot}\mathcal{C}\vert$ is easier to handle.

In the third section, following Stasheff's approach \cite[Appendices $B$ and $C$]{Bo}, we work out a detailed proof of a classification theorem for bundles with structures in a topological groupoid. 
\begin{theorem}\label{Thm2}
Given a topological groupoid $\mathcal{G}$ and a topological space $X$, there exists a $1$-$1$ correspondence between the set of homotopy classes of maps from $X$ to $\vert\operatorname{Ner}_{\cdot}\mathcal{G}^{\mathbb{N}}\vert$ and the set of homotopy classes of numerical $\mathcal{G}$-structures on $X$. 
\end{theorem} 
It is not a new theorem, and we also claim no originality for the approach presented here as it is essentially the proof of Theorem $D$, a special case of Theorem \ref{Thm2}, in \cite[Appendix $C$]{Bo}. In fact, Stasheff has indicated that his method can be applied to more general cases (see \cite[p.126]{Mo1} and \cite[Theorem $E$ in Appendix $C$]{Bo}). A similar classification theorem in terms of Milnor's construction can be found in \cite{Ha}. It is because we need a classification theorem in terms of Segal's classifying space in \cite{Wang3}, Theorem \ref{Thm2} is discussed in details here.

\subsection*{Outline of the paper}
In the second section, we review some background notions on topological groupoids, denoted by $\mathcal{G}$, and $\mathcal{G}$-structures on topological spaces. The third section discusses how Stasheff's approach \cite[Appendix $B$ and $C$]{Bo} can be generalized to arbitrary topological groupoids. The forth section, where the novelty of the paper is, is independent of the previous two sections, and we shall apply Quillen's theorem $A$ to prove Theorem \ref{Thm1} there. 
 
As the homotopy inverse $\rho:\vert\vert X_{\cdot}\vert\vert\rightarrow \vert\vert X_{\cdot}\times S_{\cdot}\vert\vert$ in \cite{tD} appears not to be well-defined, in the appendix, we use a different construction, which is due to S. Goette, to find a well-defined map $\tau:\vert\vert X_{\cdot}\vert\vert\rightarrow \vert\vert X_{\cdot}\times S_{\cdot}\vert\vert$ and explain why it is a promising candidate for a homotopy inverse to $\pi$.  
   
\subsection*{Notation and convention}
Throughout the paper, we shall use the Quillen equivalences between the model category of simplicial sets $\operatorname{sSets}$ and the Quillen model category of topological spaces $\mathpzc{Top}$ given by the singular functor and geometric realization functor:
\[\vert-\vert:\operatorname{sSets}\leftrightarrows \mathpzc{Top}:\operatorname{Sing}_{\cdot}.\]
A simplicial space is proper if and only of it is a cofibrant object in the Reedy model category associated to the Str\o n model structure on $\mathpzc{Top}$. 
 
Given $X,Y\in\mathpzc{Top}$, $[X,Y]$ denotes the set of homotopy classes of maps from $X$ to $Y$.  
 
Recall that the natural numbers $\mathbb{N}$ is an ordered set and hence can be viewed as a category.

We have chosen to work with the category of topological spaces, but the results in this paper hold for other convenient categories of topological spaces such as the category of $k$-spaces (Kelly spaces) or the category of weakly Hausdorff $k$-spaces.
\subsection*{Acknowledgment}
I would like to thank S. Goette for suggesting the approach using Quillen's theorem $A$. The construction of the map $\tau$ in Appendix is also due to him.

\section{Topological groupoids and $\mathcal{G}$-cocycles}
 
\begin{definition}\label{Topgroupoied}
A topological groupoid $\mathcal{G}$ is an internal category in $\mathpzc{Top}$ equipped with the inverse map $  
i:M(\mathcal{G})\rightarrow M(\mathcal{G})$ and the identity-assigning map $e:O(\mathcal{G})\rightarrow M(\mathcal{G})$ such that $s\circ i = t : M(\mathcal{G})\rightarrow O(\mathcal{G})$ and 
$t\circ i= s: M(\mathcal{G})\rightarrow O(\mathcal{G})$; and the diagrams below commute   
\begin{center}
\begin{tikzpicture}
\node(Lu) at (0,2) {$M(\mathcal{G})$};
\node(Ll) at (0,0) {$O(\mathcal{G})$};
\node(Mu) at (4,2) {$M(\mathcal{G})\times M(\mathcal{G})$};
 
\node(Ru) at (8,2) {$M(\mathcal{G})\times_{O(\mathcal{G})}M(\mathcal{G})$};
\node(Rl) at (8,0) {$M(\mathcal{G})$};

\path[->, font=\scriptsize,>=angle 90]

(Ll) edge node [above]{$e$} (Rl)
(Lu) edge node [above]{$\operatorname{D}$} (Mu)
(Mu) edge node [above]{$(\operatorname{id},i)$}(Ru)  
(Lu) edge node [right]{$s$}(Ll)
   
(Ru) edge node [right]{$\circ$}(Rl);

\end{tikzpicture}
\end{center}

\begin{center}
\begin{tikzpicture}
\node(Lu) at (0,2) {$M(\mathcal{G})$};
\node(Ll) at (0,0) {$O(\mathcal{G})$};
\node(Mu) at (4,2) {$M(\mathcal{G})\times M(\mathcal{G})$};
 
\node(Ru) at (8,2) {$M(\mathcal{G})\times_{O(\mathcal{G})}M(\mathcal{G})$};
\node(Rl) at (8,0) {$M(\mathcal{G})$};

\path[->, font=\scriptsize,>=angle 90]

(Ll) edge node [above]{$e$} (Rl)
(Lu) edge node [above]{$\operatorname{D}$} (Mu)
(Mu) edge node [above]{$(i,\operatorname{id})$}(Ru)  
(Lu) edge node [right]{$t$}(Ll)
   
(Ru) edge node [right]{$\circ$}(Rl);

\end{tikzpicture}
\end{center}  
where $s$ and $t$ are the source and target maps, respectively, $\operatorname{D}$ is the diagonal map $\operatorname{D}(x)=(x,x)$ and $\circ:M(\mathcal{G})\times_{O(\mathcal{G})}M(\mathcal{G})\rightarrow M(\mathcal{G})$ is the composition map.
\end{definition}


Given a topological groupoid $\mathcal{G}$ and a topological space $X$, we can define a $\mathcal{G}$-structure on $X$.
 
\begin{definition}
\begin{enumerate}
\item A $\mathcal{G}$-cocycle on $X$ is a collection $\{ \mathfrak{U}_{\alpha};f_{\beta\alpha}\}_{\alpha,\beta\in I}$, where $\{\mathfrak{U_{\alpha}}\}_{\alpha\in I}$ is an open cover of $X$ and $f_{\beta\alpha}$ is a map 
\[f_{\beta\alpha}:\mathfrak{U}_{\alpha}\cap\mathfrak{U}_{\beta}\rightarrow M(\mathcal{G})\]
that satisfies 
\[f_{\gamma\beta}\circ f_{\beta\alpha}=f_{\gamma\alpha}.\]
In particular, $f_{\alpha\alpha}$ factors through $O(\mathcal{G})$, meaning 
\[f_{\alpha\alpha}:\mathfrak{U}_{\alpha}\rightarrow O(\mathcal{G})\xrightarrow{e} M(\mathcal{G}).\]
Hence we often omit the repetition of $\alpha$ and just write $f_{\alpha}$, thinking of it as a map from $\mathfrak{U}_{\alpha}$ to $O(\mathcal{G})$.

\item Two $\mathcal{G}$-cocycles $\{\mathfrak{U}_{\alpha};f_{\beta\alpha}\}_{\alpha,\beta\in I}$ and $\{\mathfrak{V}_{\gamma};g_{\delta\gamma}\}_{\gamma,\delta\in J}$ are isomorphic if and only if there exists a map
\[\phi_{\gamma\alpha}:\mathfrak{U}_{\alpha}\cap \mathfrak{V}_{\gamma}\rightarrow M(\mathcal{G}),\]
for each $\alpha\in I$ and $\gamma\in J$, such 
that 
\[g_{\delta\gamma}\circ \phi_{\gamma\alpha}=\phi_{\delta\beta}\circ f_{\beta\alpha},\]
or equivalently, the union 
\[\{\mathfrak{U}_{\alpha},\mathfrak{V}_{\gamma} ; f_{\beta\alpha},g_{\delta\gamma},\phi_{\gamma\alpha}\}_{\alpha,\beta\in I;\gamma,\delta\in J}\] 
constitutes a $\mathcal{G}$-cocycle on $X$. Note that the index sets $I$ and $J$ are often omitted when there is no risk of confusion.

An isomorphism class of $\mathcal{G}$-cocycles on $X$ is called a $\mathcal{G}$-structure, and the set of $\mathcal{G}$-structures on $X$ is denoted by $H^{1}(X,\mathcal{G})$.

\item Two $\mathcal{G}$-structures $u,v\in H^{1}(X,\mathcal{G})$ are said to be homotopic if and only if there exists a $\mathcal{G}$-structure $w\in H^{1}(X\times I,\mathcal{G})$ such that $i_{0}^{\ast}w=u$ and $i_{1}^{\ast}w=v$, where 
$i_{0}:X=X\times\{0\}\hookrightarrow X\times I$ and 
$i_{1}:X=X\times\{1\}\hookrightarrow X\times I$. 
$\mathcal{G}(X)$ denotes the set of homotopy classes of $\mathcal{G}$-structures, and it is a contravariant functor from $\mathpzc{Top}$ to $\operatorname{Sets}$, the category of sets.
\end{enumerate}
\end{definition}

\begin{remark}\label{isohomequivalencerelation}
In this remark, we shall expand on the definition above.    
\begin{enumerate}
\item 
\emph{Isomorphisms of $\mathcal{G}$-cocycles constitute an equivalence relation on the set of $\mathcal{G}$-cocycles.} Suppose the $\mathcal{G}$-cocycles $\{\mathfrak{U}_{\alpha};f_{\beta\alpha}\}$ and $\{\mathfrak{V}_{\gamma};g_{\delta\gamma}\}$ are isomorphic through $\phi_{\gamma\alpha}$ and the $\mathcal{G}$-cocycles $\{\mathfrak{V}_{\gamma};g_{\delta\gamma}\}$ and $\{\mathfrak{W}_{\epsilon};h_{\eta\epsilon}\}$ are isomorphic through $\psi_{\gamma\epsilon}$---namely,   
\[\{\mathfrak{U}_{\alpha};f_{\beta\alpha}\}\overset{\phi_{\gamma\alpha}}{\simeq}\{\mathfrak{V}_{\gamma};g_{\delta\gamma}\}\overset{\psi_{\epsilon\gamma}}{\simeq}\{\mathfrak{W}_{\epsilon};h_{\eta\epsilon}\}.\]
Then we can define 
\[\rho_{\epsilon\alpha,\gamma}:= \psi_{\epsilon\gamma}\circ\phi_{\gamma\alpha}: \mathfrak{U}_{\alpha}\cap \mathfrak{W}_{\epsilon} \cap \mathfrak{V}_{\gamma}\rightarrow M(\mathcal{G}),\]
for each $\alpha, \gamma$ and $\epsilon$.
Since they are compatible when $\gamma$ varies, there is a well-defined map 
\[\rho_{\epsilon\alpha}:\mathfrak{U}_{\alpha}\cap \mathfrak{W}_{\epsilon}\rightarrow M(\mathcal{G}).\]
On the other hand, from the definition of $\rho_{\epsilon\alpha,\gamma}$, we have the identity  
\[h_{\eta\epsilon}\circ\rho_{\epsilon\alpha,\delta}=\rho_{\eta\beta,\gamma}\circ f_{\beta\alpha}\] 
on
\[\mathfrak{W}_{\eta}\cap \mathfrak{U}_{\alpha}\cap \mathfrak{W}_{\epsilon}\cap \mathfrak{U}_{\beta}\cap \mathfrak{V}_{\gamma}\cap\mathfrak{V}_{\delta},\] 
for all $\gamma$ and $\delta$, and hence $\{\rho_{\epsilon\alpha}\}$ constitute an isomorphism between the $\mathcal{G}$-cocycles $\{\mathfrak{U}_{\alpha};f_{\beta\alpha}\}$ and $\{\mathfrak{W}_{\epsilon};h_{\eta\epsilon}\}$. 
 
\item \emph{The notion of homotopy between $\mathcal{G}$-structures on $X$ gives an equivalence relation on $H^{1}(X,G)$.} To see this, it suffices to observe that, for any two isomorphic $\mathcal{G}$-cocycles \[\{\mathfrak{U}_{\alpha};f_{\beta\alpha}\}_{\alpha,\beta \in I}\overset{\phi_{\gamma\alpha}}{\simeq} \{\mathfrak{U}^{\prime}_{\gamma};f^{\prime}_{\delta\gamma}\}_{\gamma,\delta\in J},\] 
there is a $\mathcal{G}$-cocycle $\{\hat{U}_{\mu},\hat{f}_{\nu\mu}\}_{\nu,\mu\in I\cup J}$ on $X\times I$ which is given by
\begin{align*}
\hat{\mathfrak{U}}_{\alpha}&:= \mathfrak{U}_{\alpha}\times (1/3,1]\\
\hat{\mathfrak{U}}_{\gamma}&:= \mathfrak{U}^{\prime}_{\gamma}\times [0,2/3).
\end{align*}
and 
\begin{align*}
\hat{f}_{\beta\alpha}(x,t)&:=   f_{\beta\alpha}(x)& t>1/3 \text{  on  }\hat{\mathfrak{U}}_{\alpha}\cap \hat{\mathfrak{U}}_{\beta}\\
\hat{f}_{\gamma\delta}(x,t)&:=  f_{\gamma\delta}(x)& t<2/3 \text{  on  } \hat{\mathfrak{U}}_{\delta}\cap \hat{\mathfrak{U}}_{\gamma}\\
\hat{\phi}_{\gamma\alpha}(x,t)&:=   \phi_{\gamma\alpha}(x)& 1/3<t<2/3 \text{  on  }\hat{\mathfrak{U}}_{\alpha}\cap \hat{\mathfrak{U}}^{\prime}_{\gamma}.
\end{align*} 
\end{enumerate}
\end{remark} 

To define numerable $\mathcal{G}$-structures on a topological space. We first recall the definition of a partition of unity. 

\begin{definition}
Given a topological space $X$ and an open cover $\{\mathfrak{U}_{\alpha}\}_{\alpha\in I}$, then $\{\lambda_{\alpha}\}_{\alpha\in I}$ is a partition of unity subordinate to the open cover $\{\mathfrak{U}_{\alpha}\}_{\alpha\in I}$ if and only if 
\begin{enumerate}
\item $\lambda_{\alpha}$ is a map $\lambda_{\alpha}:X\rightarrow [0,1]$ with
$\operatorname{supp}(\lambda_{\alpha})\subset \mathfrak{U}_{\alpha}$, for each $\alpha\in I$.
\item For every point $x\in X$, there exists a neighborhood $\mathfrak{V}_{x}$ of $x$ such that, when restricted to this neighborhood, $\lambda_{\alpha}=0$ for all but finite $\alpha\in I$.    
\item For every $x\in X$,  
\[\sum_{\alpha\in I}\lambda_{\alpha}(x)=1.\]

\end{enumerate}

A numerable open cover is an open cover that admits a partition of unity subordinate to it.

\end{definition}

Not every open cover admits a partition of unity, for example, the line with two origins.

\begin{definition}
\begin{enumerate}
\item
A numerable $\mathcal{G}$-cocycle on $X$ is a $\mathcal{G}$-cocycle $\{\mathfrak{U}_{\alpha};f_{\beta\alpha}\}$ with $\{\mathfrak{U}_{\alpha}\}$ a numerable open cover.
\item
Two numerable $\mathcal{G}$-cocycles are isomorphic if and only if they are isomorphic as $\mathcal{G}$-cocycles. An isomorphism class of numerable $\mathcal{G}$-cocycles is called a numerable $\mathcal{G}$-structure on $X$. The set of $\mathcal{G}$-structures on $X$ is denoted by $H^{1}_{nu}(X,\mathcal{G})$. 
\item
Two numerable $\mathcal{G}$-structures are homotopic if and only if they are homotopic as $\mathcal{G}$-structures via a numerable $\mathcal{G}$-structure on $X\times I$. The set of homotopy classes of numerable $\mathcal{G}$-structures is denoted by $\mathcal{G}_{nu}(X)$, and it is a contravariant functor from $\mathpzc{Top}$ to $\operatorname{Sets}$.
\end{enumerate}

\end{definition}

The following lemmas imply that, when $\mathcal{G}$ is a topological group, $\mathcal{G}_{nu}(X)=H^{1}_{nu}(X,\mathcal{G})$. 

\begin{lemma}\label{HoimpliesIso}
Let $\mathcal{G}$ be a topological group and assume
$w$ is a numerable principle $\mathcal{G}$-bundle on $X\times I$. Then there is a bundle morphism $w\rightarrow \pi^{\ast}w$, where $\pi$ is the composition $X\times I\xrightarrow{(x,t)\mapsto (x,0)} X\times\{0\}\hookrightarrow X\times I$
\end{lemma}
\begin{proof}
See \cite[Theorem $9.6$ in Chapter $4$]{Hu}.
\end{proof}
\begin{lemma}
Let $\mathcal{G}$ be a topological group. Then there is a $1$-$1$ correspondence between $H^{1}(X,\mathcal{G})$ and $\{\text{ principle } \mathcal{G}\text{-bundles }\}/\text{iso.}$.  
\end{lemma}
\begin{proof}
See \cite[Theorem $11.16$]{Sw}.
\end{proof}

The next lemma explains why it is called ``numerable open cover". This technical lemma is very useful in simplifying proofs.
    
\begin{lemma}\label{NumerableCountablecover}
Given a numerable open cover $\{\mathfrak{V}_{j}\}$ of a topological space $X$, there exists a countable numerable open cover $\{\mathfrak{W}_{n}\}_{n\in\mathbb{N}}$ such that, for each $n$, $W_{n}$ is an union of some open sets, each of which is contained in some members of the original open cover $\{\mathfrak{V}_{j}\}$. 
\end{lemma}
\begin{proof}
This lemma is due to Milnor (see \cite[Thm.7.27]{Ja} or \cite[Proposition $12.1$]{Hu} for a detailed proof).
\end{proof} 

A $\mathcal{G}$-cocycle $\{\mathfrak{U}_{\alpha},f_{\beta\alpha}\}$ on $X$ is countable and numerable if and only if $\{\mathfrak{U}_{\alpha}\}$ has countably many members and is numerable. Two countable numerable $\mathcal{G}$-cocycles are isomorphic if and only if they are isomorphic as $\mathcal{G}$-cocycles. Two countable numerable $\mathcal{G}$-structures---isomorphism classes of countable numerable $\mathcal{G}$-cocycles---are homotopic if and only if they are homotopic through a countable numerable $\mathcal{G}$-structure on $X\times I$.

\begin{corollary}\label{CountablenumerableGcocycle} 
1. Given a numerable $\mathcal{G}$-structure on $X$, there exists a countable numerable $\mathcal{G}$-cocycle representing this $\mathcal{G}$-structure.

2. Given two homotopic numerable $\mathcal{G}$-structures, there exists a countable numerable $\mathcal{G}$-cocycle on $X\times I$ such that its restrictions to $X\times\{0\}$ and $X\times\{1\}$ represent the two given $\mathcal{G}$-structures.

In particular, the set of homotopy classes of numerable $\mathcal{G}$-structures is the same as the set of homotopy classes of countable numerable $\mathcal{G}$-structures.
\end{corollary}

On the other hand, in most cases, there is no loss of generality by assuming the open cover in a $\mathcal{G}$-cocycle is numerable. The ensuing corollary results from the fact that every open cover of a paracompact Hausdorff space admits a subordinate partition of unity.  
  
\begin{corollary}
Suppose $X$ is paracompact Hausdorff, then every $\mathcal{G}$-cocycle on $X$ is a numerable $\mathcal{G}$-cocycle, and any two numerable $\mathcal{G}$-structures on $X$ are homotopic if and only if they are homotopic through a numerable $\mathcal{G}$-structure on $X\times I$. In other words, classifying $\mathcal{G}$-cocycles on a paracompact Hausdorff space $X$ is equivalent to classifying numerable $\mathcal{G}$-cocycles on $X$. 
\end{corollary}
\begin{proof}
It is because the product of a paracompact space and a compact space is paracompact.
\end{proof}

\section{A classification theorem}
This section discusses a classification theorem for numerable $\mathcal{G}$-structures, and the method presented here is taken from \cite[Theorem $D$ in Appendix $C$]{Bo}, where Stasheff classifies Haefliger's structures and indicates that his approach can be applied to any topological groupoids. 

Recall first the construction of classifying spaces in \cite[Appendix $B$]{Bo}. 
\begin{definition}\label{Stasheffspace}
Given $\mathcal{C}$ an internal category in $\mathpzc{Top}$, the associated classifying space is defined to be the quotient space
\[B\mathcal{C}:= \coprod_{\mathclap{\alpha:[k]\hookrightarrow \mathbb{N}}}\operatorname{Ner}_{k}\mathcal{C}\times \triangle^{k}_{\alpha}/\sim,\] 
where $\alpha$ can be viewed a strictly increasing sequence of natural numbers $\{i_{0}<i_{1}<...<i_{k}\}$, $\triangle^{k}_{\alpha}$ is the triangle in the infinite triangle $\triangle^{\infty}=\vert \operatorname{Ner}_{\cdot}\mathbb{N}\vert$ with vertices $\{i_{0},i_{1},...,i_{k}\}$, and the relation $\sim$ is given by  
\begin{align*}
(f_{10},...,f_{kk-1};t_{0},,,\overbracket[.5pt][1pt]{0}^{j},,,t_{k})&\sim (f_{10},...,f_{j+1j}\circ f_{jj-1},...,f_{kk-1};t_{0},,,t_{j-1},t_{j+1},...,t_{k}),\\
(f_{10},...,f_{kk-1};0,t_{1},...,t_{k})&\sim (f_{21},...,f_{kk-1}; t_{1},...,t_{k}),\\
(f_{10},...,f_{kk-1};t_{0},...,t_{k-1},0)&\sim (f_{10},...,f_{k-1k-2}; t_{0},...,t_{k-1}).
\end{align*} 
\end{definition}

\begin{theorem}\label{Classificationthm}
Given a topological groupoid $\mathcal{G}$,
there is a $1$-$1$ correspondence:
\[[X,B\mathcal{G}]\leftrightarrow  \mathcal{G}_{nu}(X),\]
for every $X\in \mathpzc{Top}$.
\end{theorem}

\begin{proof}
Firstly, we note there is a canonical open cover of $B\mathcal{G}$ given by the preimage $\mathfrak{U}_{j}:=t_{j}^{-1}((0,1])$, where $t_{j}$ is induced by the projection
\[\coprod_{\mathclap{\mathclap{i_{0}<...<i_{k}}}}\operatorname{Ner}_{k}\mathcal{G}\times \triangle^{k}_{i_{0}<...<i_{k}}\xrightarrow{t_{i_{s}}}[0,1].\]
The collection $\{t_{j}\}_{i\in\mathbb{N}}$ is not locally finite and hence not a partition of unity.  To construct a partition of unity with respect to $\{\mathfrak{U}_{j}\}$, we consider maps $w_{i}; v_{i}:\mathfrak{U}_{i}\times [0,1]\rightarrow [0,1]$ defined by 
\[w_{i}(x,s):=\max\{0, t_{i}(x)-s\sum_{j< i}t_{j}(x)\}\] 
and  
\begin{equation}\label{homotopyconnectingv1andt}
v_{i}(x,s):=\frac{w_{i}(x,s)}{\sum_{j=0}^{\infty}w_{j}(x,s)},
\end{equation} 
respectively. Observe that $v_{i}(0,x)=t_{i}$ and $\{v_{i}(1,x)\}$ is locally finite and constitutes a partition of unity subordinate to the open cover $\{\mathfrak{U}_{j}\}$. The universal $\mathcal{G}$-cocycle on $B\mathcal{G}$ is then given as follows:
\begin{align*}
\gamma_{i_{l}i_{j}}:\coprod_{\mathclap{i_{0}<....<i_{k}}}\operatorname{Ner}_{k}\mathcal{G}\times \triangle^{k}_{i_{0}<...<i_{k}}& \mapsto M(\mathcal{G})\\
(g_{i_{1}i_{0}},...,g_{i_{k}i_{k-1}};t_{0},...,t_{k})& \mapsto \begin{cases} g_{i_{l}i_{l-1}}\circ...\circ g_{i_{j+1}i_{j}}  & \text{ for }  i_{j}<i_{l}\\
 s(g_{i_{j+1}i_{j}}) & \text{ for } i_{j}=i_{l}\\
  (g_{i_{l}i_{l-1}}\circ...\circ g_{i_{j+1}i_{j}})^{-1} & \text{ for } i_{j}>i_{l} 
  \end{cases} 
\end{align*}
 
On the other hand, given a countable numerable open cover $\{\mathfrak{U}_{\alpha}\}$ of a topological space $X$ and a subordinate partition of unity $\{\lambda_{\alpha}\}$, one can define the space 
\[X_{\mathfrak{U}}:=\coprod_{\mathclap{\alpha:[k]\rightarrow \mathbb{N}}}\mathfrak{U}_{\alpha}\times \triangle^{k}_{\alpha}/\sim,\]
where $\alpha=\{i_{0}<...<i_{k}\}\subset \mathbb{N}$
and the equivalence relation is
\[(x; i_{0},...,\overbracket[.5pt][1pt]{0}^{j},...,i_{k})\sim(x; i_{0},...,i_{j-1},i_{j+1},...,i_{k}),\]
for any $x\in \mathfrak{U}_{i_{0}...i_{j}...i_{k}}$.
There is a homotopy equivalence $\lambda$ from $X$ to $X_{\mathfrak{U}}$ given by
\begin{align*}
\lambda:X &\rightarrow X_{\mathfrak{U}}\\
x\in \mathfrak{U}_{i_{0}...i_{k}}&\mapsto (x;\lambda_{i_{0}}(x),....,\lambda_{i_{k}}(x)),
\end{align*}
whose homotopy inverse is the canonical projection  
\[p:X_{\mathfrak{U}}\rightarrow X.\] 
It is clear that $p\circ\lambda=\operatorname{id}_{X}$ and there is an obvious linear homotopy connecting $\operatorname{id}_{X_{\mathfrak{U}}}$ and $\lambda\circ p$. In this way, we see that the homotopy type of $X_{\mathfrak{U}}$ is independent of partition of unities.

Now, in view of Corollary \ref{CountablenumerableGcocycle}, we may assume all numerable $\mathcal{G}$-cocycles on $X$ are countable, and they are homotopic if they are homotopic through a countable numerable $\mathcal{G}$-cocycles on $X\times I$. Given a countable numerable $\mathcal{G}$-cocycle $u=\{\mathfrak{U}_{\alpha},g_{\beta\alpha},\lambda_{\alpha}\}$, we have the composition
\begin{equation}\label{InducedmapinClassificationThm}
X\xrightarrow{\lambda} X_{\mathfrak{U}}\xrightarrow{Bu} B\mathcal{G},
\end{equation}   
where the map $Bu$ is induced from the assignment
\[(x;t_{i_{0}},...,t_{i_{k}})\mapsto (g_{i_{1}i_{0}}(x),g_{i_{2}i_{1}}(x),...,g_{i_{k}i_{k-1}}(x);t_{i_{0}},t_{i_{1}},...,t_{i_{k}}).\]
Suppose two $\mathcal{G}$-cocycles $u=\{\mathfrak{U}_{\alpha},g_{\beta\alpha},\lambda_{\alpha}\}$ and $v=\{\mathfrak{V}_{\gamma},f_{\delta\gamma},\mu_{\gamma}\}$ on $X$ are homotopic, then their induced maps $Bu\circ\lambda$ and $Bv\circ\mu$ are also homotopic. This can be seen from the diagram below: 
\begin{center} 
\begin{tikzpicture}
\node(Lu) at (0,4) {$X$};
\node(Lm) at (0,2) {$X\times I$};
\node(Ll) at (0,0) {$X$}; 
\node(Ru) at (3,4) {$X_{\mathfrak{U}}$};
\node(Rm) at (3,2) {$X_{\mathfrak{W}}$};
\node(Rl) at (3,0) {$X_{\mathfrak{V}}$};
\node(RRm) at (6,2){$B\mathcal{G}$};

\path[->, font=\scriptsize,>=angle 90]

(Lu) edge node [above]{$\lambda$}(Ru)  
(Lu) edge node [right]{$\iota_{0}$}(Lm)
(Ll) edge node [right]{$\iota_{1}$}(Lm)
(Lm) edge node [above]{$\nu$}(Rm) 
(Ll) edge node [above]{$\mu$}(Rl) 
(Ru) edge (Rm)
(Rm) edge (Rl)
(Ru) edge node [yshift=.5em,xshift=.5em]{$Bu$}(RRm)
(Rm) edge node [above]{$Bw$}(RRm)
(Rl) edge node [yshift=.5em,xshift=-.5em]{$Bv$}(RRm);
\end{tikzpicture}
\end{center}
where $w=\{\mathfrak{W}_{\epsilon},h_{\eta\epsilon},\nu_{\epsilon}\}$ is a countably numerable $\mathcal{G}$-cocycle on $X\times I$ connecting $u$ and $v$, meaning $\iota^{\ast}_{0}w=u$ and $\iota^{\ast}_{1}w=v$ (see Remark \ref{isohomequivalencerelation}). We may also assume that $\{\iota^{\ast}_{0}\nu_{\epsilon}\}=\{\lambda_{\alpha}\}$ and $\{\iota^{\ast}_{1}\nu_{\epsilon}\}=\{\mu_{\alpha}\}$.
   
Thus, there is a well-defined map of sets 
\[\Psi:\mathcal{G}_{nu}(X)\rightarrow [X,B\mathcal{G}].\]
Since $\mathcal{G}_{nu}(X)$ is a contravariant functor, given any map $X\rightarrow B\mathcal{G}$, by pulling back the universal $\mathcal{G}$-cocycle on $B\mathcal{G}$, one obtains a $\mathcal{G}$-cocycle on $X$. It is also clear that pullback $\mathcal{G}$-cocycles along homotopic maps are homotopic, and hence there is a well-defined map of sets
\[\Phi:[X,B\mathcal{G}]\rightarrow\mathcal{G}_{nu}(X).\]

To see $\Phi$ is the inverse of $\Psi$, we first note that $\Psi\circ\Phi=\operatorname{id}$ is obvious as the collection $\{v_{i}(x,s)\}$ defined in equation \eqref{homotopyconnectingv1andt} connects $\{t_{i}\}$ and $\{v_{i}(1,-)\}$ and hence gives the homotopy between \[\Psi\circ\Phi(f)(x)=(g_{i_{1}i_{0}}(x),...,g_{i_{k}i_{k-1}}(x);v_{i_{0}}(1,x),...,v_{i_{k}}(1,x))\] 
and 
\[f(x)=(g_{i_{1}i_{0}}(x),...,g_{i_{k}i_{k-1}}(x);t_{i_{0}}(x),...,t_{i_{k}}(x)),\] 
for any map $f:X\rightarrow B\mathcal{G}$. Secondly, we observe that, given $u=\{\mathfrak{U}_{\alpha},g_{\beta\alpha},\lambda_{\alpha}\}$ a $\mathcal{G}$-cocycle on $X$, the pullback $\mathcal{G}$-cocycle along $Bu\circ\lambda$ is simply a restriction of $u$, namely that the pullback open cover is a subcover of $\{\mathfrak{U}_{\alpha}\}$, and thus $u$ and $(Bu\circ\lambda)^{\ast}u$ are isomorphic.   
\end{proof}

\section{Fat realization and Segal's classifying space}
In this section, we shall employ a variant of Quillen's theorem $A$ \cite[Section $1.4$]{Wa3} to show Segal's classifying space and the fat realization of an internal category in $\mathpzc{Top}$ are homotopy equivalent under some mild conditions.
  
We first recall the construction of Segal's classifying space \cite[Section $3$]{Se1}.
   
\begin{definition}
Let $\mathcal{C}$ be an internal category in $\mathpzc{Top}$. Then the associated unraveled category $\mathcal{C}^{\mathbb{N}}$ is defined as the subcategory of $\mathcal{C}\times \mathbb{N}$ obtained by deleting those morphisms $(f,i\leq i)$ when $f\neq \operatorname{id}$; and Segal's classifying space of $\mathcal{C}$ is the geometric realization of the associated unraveled category $\mathcal{C}^{\mathbb{N}}$, namely $\vert\operatorname{Ner}_{\cdot}\mathcal{C}^{\mathbb{N}}\vert$.
\end{definition} 
 
\begin{lemma}\label{Unraveledrealization}
Let $S_{\cdot}$ denote the semi-simplicial space given by
\[S_{k}:=\{\text{strictly increasing maps from $[k]$ to $\mathbb{N}$}\}.\] 
Then there are homeomorphisms
\[\vert\vert \operatorname{Ner}_{\cdot}\mathcal{C}\times S_{\cdot}\vert\vert \simeq B\mathcal{C} \simeq \vert \operatorname{Ner}_{\cdot}\mathcal{C}^{\mathbb{N}}\vert .\]
\end{lemma} 
\begin{proof}
The first homeomorphism (left) is clear as $S_{\cdot}$ is just another way of interpreting triangles $\{\triangle^{k}_{\alpha}\}$, where $\alpha$ is a strictly increasing sequence of natural numbers of length $k+1$.

For the second homeomorphism, we observe that the inclusion
\begin{align*}
\coprod_{\mathclap{\alpha:[k]\hookrightarrow \mathbb{N}}}\operatorname{Ner}_{k}\mathcal{C}\times \triangle^{k}_{\alpha} &\hookrightarrow \coprod_{k}\operatorname{Ner}_{k}\mathcal{C}^{\mathbb{N}}\times \triangle^{k}\\
(c_{0}\rightarrow...\rightarrow c_{k},t,i_{0}<...<i_{k})&\mapsto \big((c_{0},i_{0})\rightarrow...\rightarrow (c_{k},i_{k}),t\big)  
\end{align*}
descends to a homeomorphism
\[B\mathcal{C}\xrightarrow{\approx} \vert \operatorname{Ner}_{\cdot}\mathcal{C}^{\mathbb{N}}\vert.\]
\end{proof}  

Now we can state our main theorem (compare with \cite[Proposition $2$]{tD}).  
\begin{theorem}\label{tomDieckthm}
The canonical projection 
\[B\mathcal{C}\simeq \vert\operatorname{Ner}_{\cdot}\mathcal{C}^{\mathbb{N}}\vert\simeq \vert\vert\operatorname{Ner}_{\cdot}\mathcal{C}\times S_{\cdot}\vert\vert\xrightarrow{\pi} \vert\vert \operatorname{Ner}_{\cdot}\mathcal{C}\vert\vert\]
is a homotopy equivalence, provided $\operatorname{Ner}_{k}\mathcal{C}$ has the homotopy type of a $\operatorname{CW}$-complex, for every $k$. 
\end{theorem}   
\begin{proof}  
Firstly, we observe that there is a commutative diagram for any simplicial space $X_{\cdot}$ with $X_{k}$ having the homotopy type of a $\operatorname{CW}$-complex, for each $k$:

\begin{center}
\begin{equation}\label{Diagramofheqs}
\begin{tikzpicture}[baseline=(current  bounding  box.center)]

\node(Lu) at (2,3) {$\vert X_{\cdot}^{\mathbb{N}}\vert$};
\node(Lm) at (0,2) {$\vert X_{\cdot}^{\mathbb{N},p} \vert$}; 
\node(Ru) at (7.5,3) {$\vert\vert X_{\cdot}\times S_{\cdot}\vert\vert$};
\node(Rm) at (5,2) {$\vert\vert (X_{\cdot}\times S_{\cdot})^{p} \vert\vert$};

\node(Rmm) at (7.5,1){$\vert\vert X_{\cdot}\vert\vert$};

\node(Rl) at (5,0) {$\vert\vert  X_{\cdot}^{p} \vert\vert$};
 
\node(Rll)at (5,-2){$\vert  X_{\cdot}^{p} \vert$};

\path[->, font=\scriptsize,>=angle 90]

 (Lm) edge node [yshift=.3em]{$\sim$}(Lu)
 (Rm) edge node [yshift=.3em]{$\sim$}(Ru)
 (Rm) edge node [right]{$\beta$}(Rl)
 (Rl) edge node [yshift=.3em]{$\sim$}(Rmm)
 (Rl) edge node [right]{$\wr$}(Rll)
 (Ru) edge node [right]{$\pi$}(Rmm);

\draw [->](Lm) to [out=-90,in=180] node [yshift=.-1em]{\scriptsize $\gamma$}(Rll); 
\draw [double equal sign distance](Lu) to (Ru);  
\draw [double equal sign distance](Lm) to (Rm);

\end{tikzpicture}
\end{equation}
\end{center}
where arrows with the symbol $\sim$ stand for homotopy equivalences and the (semi-) simplicial space $Y_{\cdot}^{p}$ is the properization of a (semi-) simplicial space $Y_{\cdot}$, namely level-wisely applying the singular functor and geometric realization to $Y_{\cdot}$:
\[Y^{p}_{k}:=\vert\operatorname{Sing}_{\cdot}Y_{k}\vert.\]
The simplicial space $Y_{\cdot}^{\mathbb{N}}$ in the diagram above is given by
\[Y_{n}^{\mathbb{N}}:=\coprod_{\mathclap{k_{0}\leq...\leq k_{n}}}Y_{l},\]
where $l$ is the number of the distinct members in $\{k_{0},...,k_{n}\}$. The degenerate map $s_{i}^{\mathbb{N}}:Y_{n-1}^{\mathbb{N}}\rightarrow Y_{n}^{\mathbb{N}}$ is given by identities, sending the copy of $Y_{l}$ indexed by $k_{0}\leq...\leq k_{i}\leq...\leq k_{n}$ to another copy indexed by $k_{0}\leq...\leq k_{i}=k_{i}\leq...\leq k_{n}$. To define its face maps, we first group together the members in $\{k_{0},...,k_{n}\}$ that are the same. Meaning, given a sequence of increasing sequence $k_{0}\leq...\leq k_{n}$ that contains $l$ distinct numbers, we partition it into $l$ groups: 
\begin{equation}\label{groupingnumbers}
\overbrace{...}^{1}<\overbrace{...}^{2}<...<\overbrace{...}^{l}.
\end{equation} 
Then we define 
\[d^{\mathbb{N}}_{i}:Y_{n}^{\mathbb{N}}\rightarrow Y_{n-1}^{\mathbb{N}}\] 
to be 
\[d^{\mathbb{N}}_{i}\vert_{Y_{l}}:=\operatorname{id}\] 
when $k_{i}$ belongs to the group of more than one member, or otherwise 
\[d^{\mathbb{N}}_{i}\vert_{Y_{l}}:= d_{j}:Y_{l}\rightarrow Y_{l-1},\] 
where $Y_{l}$ is indexed by the given sequence and $k_{i}$ belongs to the $j$-the group in figure \eqref{groupingnumbers}.
It is not difficult to see from the construction that $Y_{\cdot}\times S_{\cdot}$ can be obtained by throwing away the degenerate part of $Y_{\cdot}^{\mathbb{N}}$---namely, those components indexed by non-strictly increasing numbers. Furthermore, if $Y_{\cdot}=\operatorname{Ner}_{\cdot}\mathcal{C}$, we have $Y_{\cdot}^{\mathbb{N}}=\operatorname{Ner}_{\cdot}\mathcal{C}^{\mathbb{N}}$.  

Now, we should expand on the homotopy equivalences in diagram \eqref{Diagramofheqs}. Firstly, since $X^{\mathbb{N}}_{\cdot}$ is proper (with respect to the Str\o m model structure) and both $X_{k}^{\mathbb{N}}$ and $X_{k}^{\mathbb{N},p}$ have the homotopy type of $\operatorname{CW}$-complexes, the upper right arrow is a homotopy equivalence \cite[VII, Proposition $3.6$]{GJ}, and therefore, in view of the homeomorphism $\vert Y^{\mathbb{N}}\vert\simeq \vert\vert Y_{\cdot}\times S_{\cdot}\vert\vert$ for any simplicial space $Y_{\cdot}$, we immediately get the upper left arrow is also a homotopy equivalence. Secondly, following from the fact that $X^{p}_{\cdot}\rightarrow X_{\cdot}$ is a level-wise homotopy equivalence and Proposition $A.1$ in \cite{Se3}, we have the map $\vert\vert X_{\cdot}^{p}\vert\vert\rightarrow \vert\vert X_{\cdot}\vert\vert$ is also a homotopy equivalence. Hence, in view of diagram \eqref{Diagramofheqs}, it suffices to show the map 
\begin{equation}\label{gammamap}
\gamma:\vert X_{\cdot}^{\mathbb{N},p}\vert \rightarrow \vert X_{\cdot}^{p}\vert
\end{equation} 
is a homotopy equivalence. Because both simplicial spaces involved in the map \eqref{gammamap} are proper, it is enough to prove the map
\[\operatorname{Sing}_{k} X_{\cdot}^{\mathbb{N}}\rightarrow \operatorname{Sing}_{k}X_{\cdot}\]
induces a homotopy equivalence, for every $k$. Now, in the case where $X_{\cdot}=\operatorname{Ner}_{\cdot}\mathcal{C}$, we have $X^{\mathbb{N}}_{\cdot}$ is the nerve of $\mathcal{C}^{\mathbb{N}}$ and the map \eqref{gammamap} is given by the natural forgetful functor 
\begin{align*}
\mathcal{C}^{\mathbb{N}}\rightarrow &\mathcal{C}\\
(c,k)\mapsto & c\\
(c\rightarrow d,k\leq l)\mapsto &  c\rightarrow d
\end{align*} 
Because the nerve ($\operatorname{Ner}_{\cdot}$) and unraveling ($\mathcal{C}\mapsto\mathcal{C}^{\mathbb{N}}$) constructions commute with the singular functor, the map
\[\operatorname{Sing}_{k}\operatorname{Ner}_{\cdot}\mathcal{C} ^{\mathbb{N}}\rightarrow \operatorname{Sing}_{k}\operatorname{Ner}_{\cdot}\mathcal{C}\] 
is identical to 
\[\operatorname{Ner}_{\cdot}\operatorname{Sing}_{k}\mathcal{C}^{\mathbb{N}}\rightarrow \operatorname{Ner}_{\cdot}\operatorname{Sing}_{k}\mathcal{C}.\] 
Therefore, if one can show the functor  
\[\mathcal{C}^{\mathbb{N}}\rightarrow \mathcal{C}\]
induces a homotopy equivalence, for any discrete category $\mathcal{C}$, then we are done.

Let's pause for a moment and recall the variant of Quillen's theorem $A$ in \cite[Sec.1.4]{Wa3}: Given a map of simplicial space $f:X_{\cdot}\rightarrow Y_{\cdot}$, if, for any $y:\triangle^{n}_{\cdot}\rightarrow Y_{\cdot}$, the space $\vert f/(\triangle^{n}_{\cdot},y)_{\cdot}\vert $ is contractible, then $f$ induces a homotopy equivalence, where $f/(\triangle^{n}_{\cdot},y)_{\cdot}$ is the pullback of 
\[\triangle^{n}_{\cdot}\rightarrow Y_{\cdot}\leftarrow X_{\cdot}.\]

Using this version of Quillen's theorem $A$, we know if one can prove the space 
\[\vert \gamma/(\triangle^{n}_{\cdot},y)_{\cdot}\vert\]
is contractible, for every $y:\triangle^{n}_{\cdot}\rightarrow \operatorname{Ner}_{\cdot}\mathcal{C}$, then the theorem follows.

To show this, we note first that every simplex $y:\triangle^{n}_{\cdot}\rightarrow \operatorname{Ner}_{\cdot}\mathcal{C}$ factors through a non-degenerate one as illustrated below:
\begin{center}  
\begin{equation}\label{Nondegeneratetriangle}
\begin{tikzpicture}[baseline=(current bounding box.center)]
\node(Lu) at (0,2) {$\gamma/(\triangle^{n}_{\cdot},y=x\circ p)_{\cdot}$};
\node(Ll) at (0,0) {$\triangle^{n}_{\cdot}$};
\node(Mu) at (4,2) {$\gamma/(\triangle^{m}_{\cdot},x)_{\cdot}$};
\node(Ml) at (4,0) {$\triangle^{m}_{\cdot}$};
\node(Ru) at (8,2) {$\operatorname{Ner}_{\cdot}\mathcal{C}^{\mathbb{N}}$};
\node(Rl) at (8,0) {$\operatorname{Ner}_{\cdot}\mathcal{C}$};

\path[->, font=\scriptsize,>=angle 90]

(Ll) edge node [above]{$p$}(Ml)
(Ml) edge node [above]{$x$,non-deg.}(Rl)
(Lu) edge (Mu)
(Mu) edge (Ru)  
(Lu) edge (Ll)
(Mu) edge (Ml)  
(Ru) edge node [right]{$\gamma$}(Rl);
\draw[->] (Ll) to[out=-10,in=-170] node [below]{\scriptsize $y$} (Rl); 
\end{tikzpicture}
\end{equation}
\end{center}  
In view of commutative diagram \eqref{Nondegeneratetriangle} and the fact that $p$ is a trivial fibration and the category of simplicial sets $\operatorname{sSets}$ is a proper model category,  , we may assume $y$ is non-degenerate. In this case, $y:\triangle^{n}_{\cdot}\rightarrow \operatorname{Ner}_{\cdot}\mathcal{C}$ is induced from a functor $[n]\rightarrow \mathcal{C}$, and hence, the pullback simplicial set $\gamma/(\triangle^{n}_{\cdot},y)_{\cdot}$ can be identified with $\operatorname{Ner}_{\cdot}([n]^{\mathbb{N}})$ and the map $\gamma/(\triangle^{n}_{\cdot},y)_{\cdot}\rightarrow \triangle^{n}_{\cdot}$ can be realized by the forgetful functor $[n]^{\mathbb{N}}\xrightarrow{\pi_{0}} [n]$.

Now, we claim the forgetful functor $[n]^{\mathbb{N}}\xrightarrow{\pi_{0}} [n]$ induces a homotopy equivalence. Consider the full subcategory $[n]^{\mathbb{N},\prime}$ of $[n]^{\mathbb{N}}$ which consists of objects $(k,l)$ with $k\leq l$. There is a natural projection 
\begin{align*}
\pi_{1}:[n]^{\mathbb{N}}&\rightarrow [n]^{\mathbb{N},\prime},\\
            (k,l)&\mapsto (k,k)& k\geq l,\\
            (k,l)&\mapsto (k,l)& k\leq l.            
\end{align*}
Suppose $\iota:[n]^{\mathbb{N},\prime}\rightarrow [n]^{\mathbb{N}}$ is the canonical inclusion, then $\pi_{1}\circ\iota=\operatorname{id}$ is obvious, and on the other hand, there is a natural transformation $\phi_{1}:\operatorname{id}\mapsto \iota\circ\pi_{1}$ given by
\begin{align*}
\phi_{1}(k,l):(k,l)&\rightarrow (k,k)& k\geq l\\
              (k,l)&\xrightarrow{\operatorname{id}} (k,l) & k\leq l.  
\end{align*}  
Thus, the functors $\iota$ and $\pi_{1}$ are inverse equivalences of categories.
 
Similarly, there is a natural projection 
\begin{align*}
\pi_{2}:[n]^{\mathbb{N},\prime}&\rightarrow [n]\\
                            (k,l)&\mapsto k,
\end{align*}
which, along with the canonical inclusion 
\begin{align*}
\iota:[n]&\rightarrow [n]^{\mathbb{N},\prime}\\
k&\mapsto (k,k),\\
\end{align*}  
gives an equivalence of categories. More precisely, we have $\pi_{2}\circ\iota=\operatorname{id}$ and the natural transformation $\phi_{2}:\iota\circ\pi_{2}\mapsto \operatorname{id}$ given by  
\[\phi_{2}(k,l):(k,k)\rightarrow (k,l),\]
for every $(k,l)\in [n]^{\mathbb{N},\prime}$.
 
Thus, we have shown the commutative diagram of equivalences of categories:
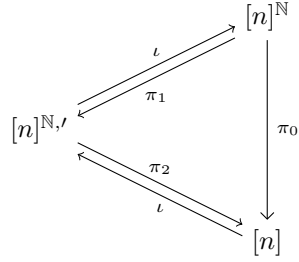
\begin{center}
\begin{tikzpicture}
\node(Lm) at (0,1.5) {$[n]^{\mathbb{N},\prime}$};  
\node(Ru) at (3,3) {$[n]^{\mathbb{N}}$};
\node(Rl) at (3,0) {$[n]$};

\path[->, font=\scriptsize,>=angle 90] 

(Ru) edge node [right]{$\pi_{0}$}(Rl);

\draw[transform canvas={yshift=0.5ex},->] (Lm) --(Ru) node[above,midway] {\scriptsize $\iota$};
\draw[transform canvas={yshift=-0.5ex},->](Ru) -- (Lm) node[below,midway,yshift=-.2em] {\scriptsize $\pi_{1}$};

\draw[transform canvas={yshift=0.5ex},->] (Lm) --(Rl) node[above,midway] {\scriptsize $\pi_{2}$};
\draw[transform canvas={yshift=-0.5ex},->](Rl) -- (Lm) node[below,midway] {\scriptsize $\iota$};

\end{tikzpicture}
\end{center}
and in particular, the space $\vert \operatorname{Ner}_{\cdot}[n]^{\mathbb{N}}\vert$ is contractible. 
 
\end{proof}
 
\addtocontents{toc}{\protect\setcounter{tocdepth}{0}} 
\begin{appendices}
\renewcommand\thesection{}
   
\section*{A candidate for a homotopy inverse to $\pi$}
\addtocontents{toc}{\protect\setcounter{tocdepth}{1}}
\addcontentsline{toc}{section}{\hspace*{2.2em} A candidate for a homotopy inverse to $\pi$}

In \cite[p.7]{tD}, in order to prove the projection  
\[\pi:\vert\vert X_{\cdot}\times S_{\cdot}\vert\vert\rightarrow \vert\vert X_{\cdot}\vert\vert\]
is a homotopy equivalence, a map $\rho$ is constructed, and it is given by the assignment 
\begin{multline}\label{tomDieckmap}
(y;t_{0},...,t_{n})\in X_{n}\times \triangle^{n} \\
\mapsto (y,1<...<n;s_{1,n}(t_{0},...,t_{n}),...,s_{n,n}(t_{0},...,t_{n}))\in X_{n}\times S_{n}\times \triangle^{n}, 
\end{multline}  
where 
\[s_{j,n}(t_{0},...,t_{n}):=(j+1)\sum_{E}\max(0,\min_{j\in E}t_{j}-\max_{j\not\in E}t_{j})\]
and $E$ runs through all subsets of $[n]$ with $j+1$ elements. However, this assignment does not respect face maps. In fact, the first and second components in the assignment \eqref{tomDieckmap} should also depend on $t_{0},...,t_{n}$.
 
Here we present a construction of a well-defined map $\vert\vert X_{\cdot}\vert\vert\rightarrow\vert\vert X_{\cdot}\times S_{\cdot}\vert\vert$ and conjecture that it is a homotopy inverse to the map $\pi$. 
  
\begin{theorem*}
There exists a map $\tau:\vert\vert X_{\cdot}\vert\vert\rightarrow\vert\vert X_{\cdot}\times S_{\cdot}\vert\vert$ such that $\pi\circ \tau$ is homotopy equivalent to $\operatorname{id}$. 
\end{theorem*}
\begin{proof}
To define the map $\tau$, we first recall that, given a simplex $\triangle^{n}:=\vert\operatorname{Ner}_{\cdot}[n]\vert$, any sequences of subsets of $[n]$
\[A_{0}\subset...\subset A_{k}\] represents a $k$-simplex in the barycentric subdivision of $\triangle^{n}$, denoted by $\operatorname{Sd}\triangle^{n}$. Then, we consider the following assignment   
\begin{align*}
\tau_{n,k}:X_{n}\times \operatorname{Sd}_{k}\triangle^{n}&\rightarrow  X_{k}\times S_{k}\times \triangle^{k}\\
(x,A_{0}\subset...\subset A_{k},t)&\mapsto (u^{\ast}x,\vert A_{0}\vert < \vert A_{1}\vert <...<\vert A_{k}\vert,t), 
\end{align*} 
where $\operatorname{Sd}_{k}\triangle^{n}$ stands for the set of $k$-simplices in $\operatorname{Sd}\triangle^{n}$, $t$ is a point in a $k$-simplex, and $\vert A_{i}\vert$ is the size of $A_{i}$. The map $u^{\ast}$ is given by the assignment
\begin{align*}
[k]&\rightarrow [n]\\
(1,2,...,k)&\mapsto(\max(A_{0}),\max(A_{1}),...,\max(A_{k}))\subset[n],
\end{align*} 
and $\max(A_{i})$ stands for the maximal element in the set $A_{i}$. It is clear that $\{\tau_{n,k}\}_{k\in [n]}$ induces a map 
\[\tau_{n}:X_{n}\times \operatorname{Sd}\triangle^{n}\rightarrow \vert\vert X_{\cdot}\times S_{\cdot}\vert\vert.\]

To see it respects the face maps in $X_{\cdot}$, we assume $y=d_{i}^{\ast}x\in X_{n-1}$ and express the image of the simplex   
\[(x,A_{0}\subset A_{1}\subset...\subset A_{k})\]
with $A_{k}\subset [n]\setminus \{i\}$ under $\tau_{n}$ in $\vert\vert X_{\cdot}\times S_{\cdot}\vert\vert$ by
\begin{equation}\label{Eq1:2nd}
(u^{\ast}x,\vert A_{0}\vert<...<\vert A_{k}\vert)
\end{equation} 
and the image of the simplex 
\[(y,B_{0},...,B_{k})\]
with $d_{i}\vert_{B_{j}}:B_{j}\xrightarrow{\simeq} A_{j}$ for each $j$ under $\tau_{n-1}$ in $\vert\vert X_{\cdot}\times S_{\cdot}\vert\vert$ by 
\begin{equation}\label{Eq2:2nd}
(v^{\ast}y,\vert B_{0}\vert<...<\vert B_{k}\vert).
\end{equation}

Now, the second components in simplices \eqref{Eq1:2nd} and \eqref{Eq2:2nd} are clearly the same, in view of the assumption 
\[d_{i}\vert_{B_{j}}:B_{j}\xrightarrow{\simeq} A_{j},\]
and the same assumption also implies the compositions
\begin{align*}
(1,2,...,k)&\mapsto (\max(B_{0}),...,\max(B_{k}))\subset [n-1]\xrightarrow{d_{i}}[n];\\
(1,2,....,k)&\mapsto (\max(A_{0}),...,\max(A_{k}))\subset [n]
\end{align*}
are identical. Hence, we have $u^{\ast}x=v^{\ast}d_{i}^{\ast}x=v^{\ast}y$, meaning that the first components in simplices \eqref{Eq1:2nd} and \eqref{Eq2:2nd} are also identical.

The homotopy between $\pi\circ \tau\simeq \operatorname{id}$ is very easy to describe. It is given by the linear homotopy from the identity to  the last vertex map. Pictorially, it looks like the following:

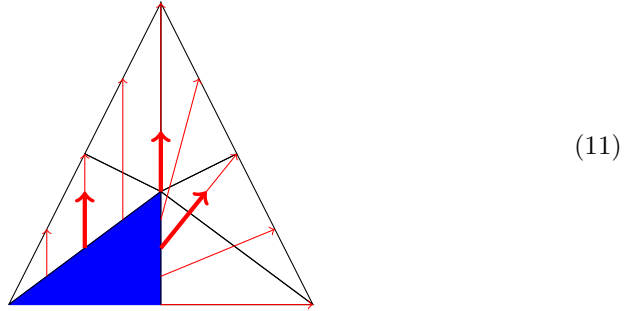
\begin{center}
\begin{equation}\label{Fig1:2nd}
\begin{tikzpicture}[baseline=(current  bounding  box.center)] 
\draw (4,0)--(3,2)--(2,1.5)--(4,0);
\draw [blue,fill=blue](0,0)--(2,0)--(2,1.5)--(0,0);
\draw (2,0)--(4,0)--(2,1.5)--(2,0);
\draw (0,0)--(2,1.5)--(1,2)--(0,0);
\draw (1,2)--(2,1.5)--(2,4)--(1,2);
\draw (3,2)--(2,1.5)--(2,4)--(3,2);

\draw [ultra thick,->,red] (1,.75) to (1,1.5);
\draw [ultra thick,->,red] (2,.75) to (2.6,1.5);
\draw [ultra thick,->,red] (2,1.5) to (2,2.3);
\draw [->, red] (2,1.5) to (2,4);
\draw [->, red] (1.5,1.125) to (1.5,3);
\draw [->, red] (1,.75) to (1,2);
\draw [->, red] (.5,.375) to (.5,1);
\draw [->, red] (2,1.125) to (2.5,3);
\draw [->, red] (2,.75) to (3,2);
\draw [->, red] (2,.375) to (3.5,1);
\draw [->, red] (2,0) to (4,0);
\end{tikzpicture}
\end{equation}
\end{center}

\end{proof} 
\begin{conjecture*}
Is $\tau$ a homotopy inverse to the map $\pi$? 
\end{conjecture*}

\begin{remark*}
Though the idea is not so complicated, we find it very hard to write down the homotopy between $\tau\circ\pi$ and $\operatorname{id}$ in details. Observe first that, given any simplicial space $Y_{\cdot}$, there is a natural filtration 
\[\emptyset\subset \vert\vert Y_{\cdot}\vert\vert_{(0)}\subset...\subset \vert\vert Y_{\cdot}\vert\vert_{(k)}\subset...\subset \vert\vert Y_{\cdot}\vert\vert\]
given by truncating the simplicial space $Y_{\cdot}$. Our strategy is to define a homotopy 
\[h_{n}:\vert\vert X_{\cdot}\times S_{\cdot}\vert\vert_{(n)}\times I\rightarrow \vert\vert X_{\cdot}\times S_{\cdot}\vert\vert_{(n+1)}\]
between $\tau\circ\pi\vert_{\vert\vert X_{\cdot}\times S_{\cdot}\vert\vert_{(n)}}$ and $\operatorname{id}$, for each $n$ such that the diagram below commutes: 
\begin{center}
\begin{equation}\label{Diag1:2nd}
\begin{tikzpicture}[baseline=(current  bounding  box.center)]  
\node(Lu) at (0,2) {$\vert\vert X_{\cdot}\times S_{\cdot}\vert\vert_{(n)}\times I$};
\node(Ll) at (0,0) {$\vert\vert X_{\cdot}\times S_{\cdot}\vert\vert_{(n+1)}\times I$}; 
\node(Ru) at (4,2) {$\vert\vert X_{\cdot}\times S_{\cdot}\vert\vert_{(n+1)}$};
\node(Rl) at (4,0) {$\vert\vert X_{\cdot}\times S_{\cdot}\vert\vert_{(n+2)}$};

\path[->, font=\scriptsize,>=angle 90]

(Lu) edge (Ru)  
(Lu) edge (Ll)
(Ll) edge (Rl) 
(Ru) edge (Rl);

\end{tikzpicture}
\end{equation}
\end{center}
Then, passing to the colimit, we get the required homotopy.

The homotopy $h_{0}$ is simply the homotopy given by the $1$-simplex $(x,0<k)$, where $x\in X_{0}$. For general $n$, we decompose the homotopy $h_{n}$ into two parts. Given a non-degenerate $x\in X_{n}$, there is a canonical embedding 
\[\triangle^{n}\times [n,+\infty)\hookrightarrow \vert\vert X_{\cdot}\times S_{\cdot}\vert\vert,\]
and the first part of $h_{n}$ is the linear homotopy to the projection 
\[\triangle^{n}\times [n+\epsilon,+\infty)\rightarrow \triangle^{n}\times \{n+\epsilon\}\] 
with any thing below $\{n+\epsilon\}$ intact, where $\epsilon>0$.
The following illustrates the case $n=1$.

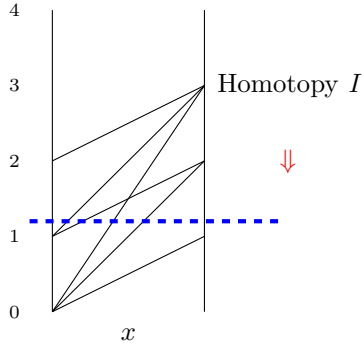
\begin{center} 
\begin{tikzpicture}
\draw (0,0) to (0,4); 
\draw (2,0) to (2,4);
\draw (0,0) to (2,1);
\draw (0,0) to (2,2);
\draw (0,0) to (2,3);
\draw (0,1) to (2,2);
\draw (0,1) to (2,3);
\draw (0,2) to (2,3);

\draw [dashed,blue,ultra thick](-.3,1.2) to (3,1.2);  
\node at (3.1,3) {$\operatorname{Homotopy}$ I}; 
\node [red] at (3.1,2) {$\Downarrow$}; 
\node at (-.5,0){\scriptsize $0$};
\node at (-.5,1){\scriptsize $1$};
\node at (-.5,2){\scriptsize $2$};
\node at (-.5,3){\scriptsize $3$};
\node at (-.5,4){\scriptsize $4$};
\node at (1,-.3){$x$};
\end{tikzpicture}
\end{center}
The second part of $h_{n}$ is nothing but a thickened version of the homotopy illustrated in figure \eqref{Fig1:2nd} except that instead of moving the vertices simultaneously, we start with the vertex of largest depth and down to the one of the smallest depth. The figure below illustrates the case $n=1$. 
\begin{center} 
\begin{tikzpicture}
\draw (0,0) to (0,1.5);
\draw (0,0) to (2,1);
\draw (2,0) to (2,1.5);
\draw[red] (0,.3)-- (.5,.3)-- (2.04,1.07)-- (2.04,.3);
\draw[red] (0,.7) -- (1.2,.7)-- (2.1,1.15)-- (2.1,.5);
 
\draw[red] (0,1.2) -- (2.15,1.2)-- (2.15,1);
   
\draw[blue,dashed] (-.3,1.25) to (2.3,1.25);
\node[red] at (1,1.5) {$\Rightarrow$};
\node[red] at (2.5,.5){$\Downarrow$}; 
\node[red] at (-.5,.5){$\Downarrow$};
\node at (-.3,0){\scriptsize $0$};
\node at (-.3,1){\scriptsize $1$};
\end{tikzpicture} 
\end{center} 

One still needs to take care of degenerate simplices in $X_{\cdot}$, and it seems to be a cumbersome task to write down the homotopy of those degenerate simplices, even though it is possible to describe it in low-dimension. We are hoping for a better way to approach this problem.  
\end{remark*}

\end{appendices}
 



\addcontentsline{toc}{section}{\hspace*{2.2em} References} 

\bibliographystyle{alpha}

\bibliography{Reference}  

\end{document}